\title{\bf{Some unlimited families of minimal surfaces of general type with the canonical map of degree 8}}
\author{
	NGUYEN BIN\\
}
\date{\today}

\newcommand{\Addresses}{{
		\bigskip
		\footnotesize
			\text{Center for Mathematical Analysis, Geometry and Dynamical Systems}\par\nopagebreak
			\text{Departamento de Matem\'{a}tica}\par\nopagebreak
			\text{Instituto Superior T\'{e}cnico}\par\nopagebreak
			\text{Universidade de Lisboa}\par\nopagebreak
			\text{Av. Rovisco Pais, 1049-001 Lisboa, Portugal.}\par\nopagebreak
		\textit{E-mail address}: \texttt{nguyenbin@tecnico.ulisboa.pt}
				
	}}

\newcommand\blfootnote[1]{%
	\begingroup
	\renewcommand\thefootnote{}\footnote{#1}%
	\addtocounter{footnote}{-1}%
	\endgroup
}

\date{}

\documentclass[10pt]{article}
\usepackage{amsmath, amsthm, amssymb, mathrsfs, amscd, amsthm, amscd,amsfonts,calligra,mathrsfs,adjustbox,lipsum}

\usepackage{indentfirst,url,xypic}
\usepackage[all]{xy}
\usepackage{url}
\usepackage{tikz}
\usepackage[colorlinks,plainpages]{hyperref}
\hypersetup{
	colorlinks=true,
	linkcolor=blue,
	filecolor=magenta,      
	urlcolor=cyan,
}

\DeclareMathOperator{\degree}{deg}
\DeclareMathOperator{\image}{im}

\newtheorem{Theorem}{Theorem}
\newtheorem{Remark}{Remark}
\newtheorem{Proposition}{Proposition }

\newtheorem{Notation}{Notation}

\newcommand{\MSC}{\textbf{Mathematics Subject Classification (2010):}}

\newcommand{\Key}{\textbf{Key words:}}

\begin{document}
\maketitle
\begin{abstract} 
	In this note, we construct nine families of projective complex minimal surfaces of general type having the canonical map of degree 8 and irregularity $ 0 $ or $ 1 $. For six of these families the canonical system has a non trivial fixed part. 
\end{abstract}

\blfootnote{\MSC{ 14J29}.}
\blfootnote{\Key{ Surface of general type, Canonical degree, Abelian cover.}}

    \section{Introduction}
     Let $ X $ be a smooth complex surface of general type (see \cite{MR1406314} or \cite{MR2030225}) and let $ \xymatrix{\varphi_{\left| K_X\right| }:X \ar@{.>}[r] & \mathbb{P}^{p_g\left( X\right)-1}} $ be the canonical map of $ X $, where $ p_g\left( X\right) = \dim\left( H^0\left( X, K_{X}\right) \right) $ is the geometric genus and $ K_{X} $ is the canonical divisor of $ X. $ A classical result of Beauville \cite[\rm Theorem 3.1]{MR553705} says that if the image of $ \varphi_{\left | K_X\right| } $ is a surface, either $ p_g\left( \image\left(\varphi_{\left | K_X\right| } \right)  \right) =0 $ or $ \image\left(\varphi_{\left | K_X\right| } \right) $ is a surface of general type. In addition, the degree $ d $ of the canonical map of $ X $ is less than or equal to 36. \\
     
     While surfaces with $ d=2 $ has been studied thoroughly by E. Horikawa in his several papers such as \cite{MR0424831}, \cite{MR0460340}, \cite{MR517773}, \cite{MR501370}, the case where $ d $ bigger than 2 remains to be one of the most interesting open problems in the theory of surfaces. Several surfaces with $ d $ bigger than 2 have been constructed, for example with $ d = 3,5,9 $ by R. Pardini \cite{MR1103913} and S.L. Tan \cite{MR1141782}, $ d = 6,8 $ by A. Beauville \cite{MR553705}, $ d=4  $ by A. Beauville \cite{MR553705}, and F.J. Gallego and B.P. Purnaprajna \cite{MR2415082}, $ d = 16 $ by U. Persson \cite{MR527234} and C. Rito \cite{MR3619737}, $ d = 12, 24 $ by C.Rito \cite{MR3663791} \cite{MR3391024}, etc. \\
     
     In the same paper \cite{MR553705}, Beauville also proved that the degree of the canonical map is less than or equal to $ 9 $ if $\chi(\mathcal O_X)\geq 31$. Later, G. Xiao showed that if the geometric genus of $ X $ is bigger than $ 132 $, the degree of the canonical map is less than or equal to $ 8 $ \cite{MR842626}. In addition, he also proved that if the degree of the canonical map is $ 8 $ and geometric genus is bigger than $ 115 $, the irregularity $ q = h^0\left( \Omega_X^1\right) $ is less than or equal to $ 3 $ (see \cite{MR904942}). Beauville constructed an unlimited family of surfaces with $ d =8 $ and arbitrarily high geometric genus \cite{MR553705}. These surfaces have irregularity $ q=3 $ and the canonical linear system of these surfaces is base point free.\\
     
     In this note, we construct nine unlimited families of surfaces with $ d = 8 $ and $ q = 0 $ or $ q = 1 $. Furthermore, for some families the canonical linear systems are not base point free. The following theorem is the main result of this note:
     \begin{Theorem}\label{the main theorem}
     	Let $ n $ be an integer number such that $ n \ge 2 $. Then there exist minimal surfaces of general type $ X $ with canonical map $ \varphi_{\left| K_X  \right|}  $ of degree 8 and the following invariants	
     	$$
     	\begin{tabular}{|c|c|c|c|}
     	\hline
     	$ K_X^2 $ &$ p_g\left( X\right)  $&$ q\left( X\right)  $ &$ \left| K_X \right| $ is base point free\\ \hline
     	$ 16n - 8 $ & $ 2n+1 $ & $ 0 $ & yes \\ \hline
     	$ 16n - 16 $ & $ 2n $ & $ 0 $ & yes \\ \hline
     	$ 16n - 16 $ & $ 2n $ & $ 1 $  & yes\\ \hline
     	$ 16n -10 $ & $ 2n $ & $ 0 $ &  no\\ \hline
     	$ 16n  $ & $ 2n+1 $ & $ 0 $ & no \\ \hline
     	$ 16n -8 $ & $ 2n $ & $ 0 $ & no \\ \hline
     	$ 16n -8 $ & $ 2n $ & $ 1 $ & no \\ \hline
     	$ 16n - 2 $ & $ 2n $ & $ 0 $ & no \\ \hline
     	$ 16n  $ & $ 2n $ & $ 1 $ & no \\ \hline
     	\end{tabular}
     	$$
     \end{Theorem}
     \noindent
     The approach to construct these surfaces is using $ \mathbb{Z}_{2}^3- $covers with some appropriate branch loci. Note that canonical maps defined by abelian covers of  $ \mathbb{P}^2 $, and in particular the abelian covers with the group $ \mathbb{Z}_2^3 $, have been studied very explicitly by Rong Du and Yun Gao in \cite{MR3217634}.

	\section{$ \mathbb{Z}_{2}^3- $coverings}
	    The construction of abelian covers was studied by R. Pardini in \cite{MR1103912}. \\
	    \noindent
	    Let $ H_{i_1,i_2,i_3} $ denote the nontrivial cyclic subgroup generated by $ \left( i_1,i_2,i_3\right) $ of $ \mathbb{Z}_{2}^3$ for all $ \left( i_1,i_2,i_3\right) \in \mathbb{Z}_{2}^3 \setminus \left( 0,0,0\right) $, and denote by $ \chi_{j_1,j_2,j_3} $ the character of $ \mathbb{Z}_{2}^3 $ defined by
	    \begin{align*}
	    \chi_{j_1,j_2,j_3}\left( a_1,a_2,a_3\right): =  e^{\left( \pi a_1j_1\right) \textbf{i}}e^{\left( \pi a_2j_2\right) \textbf{i}}e^{\left( \pi a_3j_3\right) \textbf{i}}
	    \end{align*}
	    for all $ j_1,j_2,j_3,a_1,a_2,a_3,a_4 \in \mathbb{Z}_2. $ For sake of simplicity, from now on we use notations $ D_1, D_2, D_3, D_4, D_5, D_6, D_7 $ instead of $ D_{\left( H_{0,0,1},\chi_{0,0,1}\right) }, D_{\left( H_{0,1,0},\chi_{0,1,0}\right) }$, $ D_{\left( H_{0,1,1},\chi_{0,1,0}\right) }, D_{\left( H_{1,0,0},\chi_{1,0,0}\right) }$,  $ D_{\left( H_{1,0,1},\chi_{1,0,0}\right) }, D_{\left( H_{1,1,0},\chi_{1,0,0}\right) }, D_{\left( H_{1,1,1},\chi_{1,0,0}\right) } $, respectively. For details about the building data of abelian covers and their notations, we refer the reader to Section 1 and Section 2 of R. Pardini's work (\cite{MR1103912}). From \cite[\rm Theorem 2.1]{MR1103912} we can define $ \mathbb{Z}_{2}^3- $covers as follows:
	    \begin{Proposition} \label{Construction of cover}
	    	Let $ Y $ be a smooth projective surface. Let $ L_{\chi} $ be divisors of $ Y $ such that $ L_{\chi} \not\equiv \mathcal{O}_Y $ for all nontrivial characters $ \chi \in \left( \mathbb{Z}_{2}^3\right)^{*} \setminus \left\lbrace  \chi_{0,0,0}  \right\rbrace  $. Let $ D_1, D_2, \ldots, D_7 $ be effective divisors of $ Y $ such that the branch divisor $ B:=\sum_{i=1}^{7}{D_i} $ is reduced. Then $ \left\lbrace L_{\chi}, D_j \right\rbrace_{\chi,j}$ is the building data of a $ \mathbb{Z}_{2}^3-$cover $ \xymatrix{f:X \ar[r]& Y} $ if and only if
	    	$$
	    	\begin{adjustbox}{max width=\textwidth}
	    	\begin{tabular}{l l l l l l l l }
	    	$ 2L_{1,0,0} $&$ \equiv $&$  $&$ $&$ D_{4} $&$ +D_{5 } $&$ +D_{6} $&$ +D_{7} $ \\
	    	$ 2L_{0,1,0} $&$ \equiv $&$ D_{2} $&$ +D_{3} $&$  $&$ $&$ +D_{6} $&$ +D_{7} $ \\
	    	$ 2L_{0,0,1} $&$ \equiv D_{1 } $&$ $&$ +D_{3} $&$$&$ +D_{5} $&$  $&$ +D_{7 } $ \\
	    	$ 2L_{1,1,0} $&$ \equiv $&$ D_{2 } $&$ +D_{3} $&$ +D_{4} $&$ +D_{5} $&$  $&$  $\\
	    	$ 2L_{1,0,1} $&$ \equiv D_{1} $&$ $&$ +D_{3} $&$ +D_{4} $&$  $&$ +D_{6} $&$  $ \\
	    	$ 2L_{0,1,1} $&$ \equiv D_{1} $&$ +D_{2} $&$  $&$  $&$ +D_{5} $&$ +D_{6} $&$  $ \\
	    	$ 2L_{1,1,1} $&$ \equiv D_{1} $&$ +D_{2} $&$  $&$ +D_{4} $&$ $&$ $&$ +D_{7 } $.
	    \end{tabular}
	\end{adjustbox}
	$$	
	\end{Proposition}
	
	By \cite[\rm Theorem 3.1]{MR1103912} if each $D_\sigma$ is smooth and $B $ is a simple normal crossings divisor, then the surface $X$ is smooth. \\
	
	Also from \cite[\rm Lemma 4.2, Proposition 4.2]{MR1103912} we have:
	\begin{Proposition}\label{invariants}
	Let $ \xymatrix{f: X \ar[r]& Y} $ be a smooth $  \mathbb{Z}_{2}^3- $cover with the building data $ D_1, D_2, \ldots, D_7, L_{\chi}, \forall \chi \in \left( \mathbb{Z}_{2}^3\right)^{*} \setminus \left\lbrace  \chi_{0,0,0}  \right\rbrace $. The invariants of $ X $ are as follows:
	\begin{align*}
	2K_X & \equiv f^*\left( 2K_Y + \sum\limits_{j = 1}^7 {D_j } \right) \\
	K^2_X &= 2\left( 2K_Y + \sum\limits_{j = 1}^7 {D_j } \right)^2 \\
	p_g\left( X \right) &=p_g\left( Y \right) +\sum\limits_{\chi \in \left( \mathbb{Z}_{2}^3\right)^{*} \setminus \left\lbrace  \chi_{0,0,0}  \right\rbrace }{h^0\left( L_{\chi} + K_Y \right)} \\
	\chi\left( \mathcal{O}_X \right) &= 8\chi\left( \mathcal{O}_Y \right)  +\sum\limits_{\chi \in \left( \mathbb{Z}_{2}^3\right)^{*} \setminus \left\lbrace  \chi_{0,0,0}  \right\rbrace }{\frac{1}{2}L_{\chi}\left( L_{\chi}+K_Y\right)}. 
	\end{align*}
	\end{Proposition}

	\begin{Notation} \label{the notation of a compatible point}
	We denote $ P = \left( k_1, k_2, \ldots, k_7\right)  $ when $ D_1, D_2,\ldots, D_7 $ contain $ P $ with multiplicity $ k_1, k_2, \ldots, k_7, $ respectively.
	\end{Notation}

   \section{Constructions}
   \subsection{Construction 1}
   In this section, we construct the surfaces in the first four rows of Theorem \ref{the main theorem}.
   \subsubsection{Construction and computation of invariants}
   \label{construction1.1}
   Let $ \mathbb{F}_1 $ denote the Hirzebruch surface with the negative section  $ \Delta_0 $ with self-intersection $-1$ and  let $ \Gamma $ denote a fiber of the ruling. Let $ D_2 = 2n\Gamma $ be $ 2n $ fibers in $ \mathbb{F}_1 $ and $ D_3, D_6, D_7 \in \left| 2\Delta_0 + 2\Gamma \right| $ be smooth curves in general position. Let $ \xymatrix{f: X \ar[r] & \mathbb{F}_1}  $ be a $ \mathbb{Z}^3_2- $cover with the following branch locus
   \begin{align*}
   B = D_1 + D_2 + D_3 +D_4 +D_5 + D_6 + D_7,
   \end{align*}
   \noindent
   where $ D_1 = D_4 = D_5 = 0 $. By Proposition \ref{Construction of cover}, $ L_{0,1,0} \equiv 3\Delta_0 + \left( n+3\right) \Gamma$ and $ L_{\chi} $ is equivalent to either $ 2\Delta_0 + 2 \Gamma $ or $ \Delta_0 + \left( n+1\right) \Gamma $ for all $ L_{\chi} \ne L_{0,1,0} $. Since each $ D_{\sigma} $ is smooth and $ B $ is a normal crossings divisor, $ X $ is smooth.	Moreover, by Proposition \ref{invariants}, we get
   \begin{align*}
   2K_X & \equiv f^*\left( 2\Delta_0 + 2n\Gamma \right).
   \end{align*}
   \noindent
   This implies that $ X $ is a minimal surface of general type. Furthermore, by Proposition \ref{invariants}, the invariants of $ X $ are as follows:
   \begin{align}\label{chern class 1}
   K_X^2= 8\left( 2n-1\right)  
   \end{align}
   \begin{align} \label{pg}
   p_g\left( X\right) &= h^0\left( \Delta_0 + n\Gamma\right) =2n+1
   \end{align}
   \begin{align}\label{chi}
   \chi\left( \mathcal{O}_X\right) &=2n+2.
   \end{align}
   \noindent
   From $ \left( \ref{pg}\right)  $ and $ \left( \ref{chi}\right)  $, we get $ q\left( X\right) = 0. $ \\
   
   We show that $ \left| K_X\right| $ is not composed with a pencil by considering the following double cover
   \begin{align*}
   \xymatrix{f_1: X_{1} \ar[r]& \mathbb{F}_1}
   \end{align*}
   ramifying on  $ D_2 + D_3 + D_6 + D_7 $. We have 
   \begin{align*}
   K_{X_{1}}& \equiv f_1^*\left( \Delta_0 + n\Gamma \right). 
   \end{align*}
   Because $ \left| \Delta_0 + n\Gamma \right|  $ is not composed with a pencil, $ \left| K_{X_{1}}\right|  $ is not composed with a pencil, either. This leads to the fact that $ \left| K_X\right| $ is not composed with a pencil and the degree of the canonical map is $ 8 $. Moreover, $ \degree\left( \image\varphi_{\left| K_X\right| }  \right) =2n-1 $.
   
   \subsubsection{Variations}
   Now by adding a singular point to the above branch locus, we obtain the surfaces described  in the second row of Theorem \ref{the main theorem}. In fact, by Proposition $  \ref{Construction of cover}$, a new branch locus can be formed by adding a point $ P = \left( 0, 1,1,0,0,1,1\right)  $ (see Notation \ref{the notation of a compatible point}). And we consider the $ \mathbb{Z}_2^3- $cover on $ Y $ instead of $ \mathbb{F}_1, $ where $ Y $ is the blow up of $ \mathbb{F}_1 $ at $ P $. More precisely, let $ P $ be a point in $ \mathbb{F}_1 $ such that $ D_2, D_3, D_6, D_7 $ contain $ P $ with multiplicity $ 1,1,1,1 $, respectively. Let $ Y $ be the blow up of $ \mathbb{F}_1 $ at $ P $ and $ E $ be the exceptional divisor. If we abuse notation and denote $ D_2, D_3, D_6, D_7, \Delta_0, \Gamma $ their pullbacks to $ Y $, then $ D_2 = 2n\Gamma -E  , D_3 = 2\Delta_0 + 2\Gamma  - E, D_6 = 2\Delta_0 + 2\Gamma  - E  $ and $ D_7 = 2\Delta_0 + 2\Gamma  - E $. Let $ \xymatrix{f: X \ar[r] & Y}  $ be a $ \mathbb{Z}^3_2- $cover with the following branch locus
   \begin{align*}
   B =D_1 + D_2 + D_3 + D_4 + D_5+ D_6 + D_7,
   \end{align*}
   \noindent
   where $ D_1 = D_4 = D_5 = 0 $.  The building data is as follows:
   $$	
   \begin{tabular}{l r r r}
   	$ L_{1,0,0} \equiv$ & $ 2\Delta_0 $ & $ +2\Gamma $&$  - E $\\
   	$ L_{0,1,0} \equiv$ & $3\Delta_0   $ &$ +\left( n+3\right)\Gamma$&$ - 2E $\\
   	$ L_{0,0,1} \equiv$ & $2\Delta_0   $ &$ +2\Gamma$&$ - E $\\
   	$ L_{1,1,0} \equiv$ & $ \Delta_0   $ & $ +\left( n+1\right)\Gamma $&$  - E $\\
   	$ L_{1,0,1} \equiv$ & $ 2\Delta_0   $ & $ +2\Gamma $&$  - E $\\
   	$ L_{0,1,1} \equiv$ & $ \Delta_0   $ & $ +\left( n+1\right)\Gamma $&$  - E $\\
   	$ L_{1,1,1} \equiv$ & $ \Delta_0   $ & $ +\left( n+1\right)\Gamma $&$  - E $.\\
   \end{tabular} 
   $$
   \noindent  	
   Similarly to the above, we obtain minimal surfaces of general type with
   \begin{align*}
   K^2 = 16n-16, p_g = 2n, q = 0, d= 8,
   \end{align*}
   and $ \degree\left( \image\varphi_{\left| K_X\right| }  \right) =2n-2 $. Moreover, $ \varphi_{\left| K_X\right| } $ is a morphism. \\
   
   Analogously, by Proposition $ \ref{Construction of cover},  $ a point $ \left( 0,0,0,0,0,2,2\right)  $ can be added to the original branch locus. In fact, let $ P $ be a point in $ \mathbb{F}_1 $ such that $ D_6, D_7 $ contain $ P $ with multiplicity $ 2,2 $, respectively. Let $ Y $ be the blow up of $ \mathbb{F}_1 $ at $ P $ and $ E $ be the exceptional divisor. If we abuse notation and denote $ D_2, D_3, D_6, D_7, \Delta_0, \Gamma $ their pullbacks to $ Y $, then $ D_2 = 2n\Gamma , D_3 = 2\Delta_0 + 2\Gamma , D_6 = 2\Delta_0 + 2\Gamma  - 2E  $ and $ D_7 = 2\Delta_0 + 2\Gamma  - 2E $. Let $ \xymatrix{f: X \ar[r] & Y}  $ be a $ \mathbb{Z}^3_2- $cover with the following branch locus
   \begin{align*}
   B =D_1 + D_2 + D_3 + D_4 + D_5+ D_6 + D_7,
   \end{align*}
   \noindent
   where $ D_1 = D_4 = D_5 = 0 $.  The building data is as follows:
   $$	
   \begin{tabular}{l r r r}
   	$ L_{1,0,0} \equiv$ & $ 2\Delta_0 $ & $ +2\Gamma $&$  - 2E $\\
   	$ L_{0,1,0} \equiv$ & $3\Delta_0   $ &$ +\left( n+3\right)\Gamma$&$ - 2E $\\
   	$ L_{0,0,1} \equiv$ & $2\Delta_0   $ &$ +2\Gamma$&$ - E $\\
   	$ L_{1,1,0} \equiv$ & $ \Delta_0   $ & $ +\left( n+1\right)\Gamma $&$   $\\
   	$ L_{1,0,1} \equiv$ & $ 2\Delta_0   $ & $ +2\Gamma $&$  - E $\\
   	$ L_{0,1,1} \equiv$ & $ \Delta_0   $ & $ +\left( n+1\right)\Gamma $&$  - E $\\
   	$ L_{1,1,1} \equiv$ & $ \Delta_0   $ & $ +\left( n+1\right)\Gamma $&$  - E $.\\
   \end{tabular} 
   $$
   \noindent
   We get minimal surfaces of general type with
   \begin{align*}
   K^2 = 16n-16, p_g = 2n, q = 1, d=8,
   \end{align*}
   and $ \degree\left( \image\varphi_{\left| K_X\right| }  \right) =2n-2 $. Furthermore, $ \varphi_{\left| K_X\right| } $ is a morphism. Therefore we obtain the surfaces described  in the third row of Theorem \ref{the main theorem}. The Albanese pencil of these surfaces $ \xymatrix{X \ar[r]& Alb\left( X\right) } $ is the pullback of the Albanese pencil of the intermediate surface $ Z $, where $ Z $ is obtained by the $ \mathbb{Z}_2 -$cover ramifying on $ 2L_{1,0,0} $. For details about the surfaces with $ q > 0 $, we refer the reader to the work of Mendes Lopes and Pardini \cite{MR2931875}.
   \begin{Remark}
   These surfaces in the first three rows of Theorem \ref{the main theorem} can be obtained by taking three iterated $ \mathbb{Z}_2- $covers. First, we ramify on $ D_2, D_3, D_6,$ and $ D_7 $ (i.e. $ B=2L_{0,1,0} $) and we get Horikawa's surfaces with $ K^2 = 2p_g -4 $ \cite{MR0424831}. The second cover ramifies only on nodes (i.e $ B= 2L_{1,0,0} $). These nodes come from the intersection points between $ D_2 +D_3$ and $ D_6 +D_7 $. The last cover ramifies on nodes coming from the intersection points between $ D_2 $ and $D_3$, and $ D_6 $ and $D_7 $ (i.e. $ B=2L_{0,0,1} $) (see \cite[Prop. 3.1]{MR1918136}). Moreover, the following diagram commutes
   \end{Remark}
   $$
   \xymatrix{X \ar[0,3]^{\mathbb{Z}_2^3}_f \ar[1,1]^{2L_{0,0,1}}_{f_3} \ar[3,1]_{\varphi_{\left| K_X \right| }}&&& Y\\
   	&X_2 \ar[r]^{2L_{1,0,0}}_{f_2} &X_1 \ar[2,-1]^{\varphi_{\left| K_{X_1} \right| }} \ar[-1,1]^{f_1}_{2L_{0,1,0}}&\\
   	&&&\\
   	&\mathbb{P}^{p_g\left( X\right) -1}&&} 	
   $$

   \noindent
   Now, by Proposition $ \ref{Construction of cover},  $ a point $ \left( 0,0,1,0,-1,1,2\right)  $ can be imposed on the original branch locus, where $ -1 $ in the fifth component means the exceptional divisor is added to $ D_5 $. In fact, let $ P $ be a point in $ \mathbb{F}_1 $ such that $ D_3, D_6, D_7 $ contain $ P $ with multiplicity $ 1,1,2 $, respectively. Let $ Y $ be the blow up of $ \mathbb{F}_1 $ at $ P $ and $ E $ be the exceptional divisor. If we abuse notation and denote $ D_2, D_3, D_6, D_7, \Delta_0, \Gamma $ their pullbacks to $ Y $, then $ D_2 = 2n\Gamma , D_3 = 2\Delta_0 + 2\Gamma - E, D_6 = 2\Delta_0 + 2\Gamma  - E  $ and $ D_7 = 2\Delta_0 + 2\Gamma  - 2E $. Let $ \xymatrix{f: X \ar[r] & Y}  $ be a $ \mathbb{Z}^3_2- $cover with the following branch locus
   \begin{align*}
   B =D_1 + D_2 + D_3 + D_4 + D_5+ D_6 + D_7,
   \end{align*}
   \noindent
   where $ D_1 = D_4 = 0 $ and $ D_5 = E $.  The building data is as follows:
   $$	
   \begin{tabular}{l r r r}
   	$ L_{1,0,0} \equiv$ & $ 2\Delta_0 $ & $ +2\Gamma $&$  - E $\\
   	$ L_{0,1,0} \equiv$ & $3\Delta_0   $ &$ +\left( n+3\right)\Gamma$&$ - 2E $\\
   	$ L_{0,0,1} \equiv$ & $2\Delta_0   $ &$ +2\Gamma$&$ - E $\\
   	$ L_{1,1,0} \equiv$ & $ \Delta_0   $ & $ +\left( n+1\right)\Gamma $&$   $\\
   	$ L_{1,0,1} \equiv$ & $ 2\Delta_0   $ & $ +2\Gamma $&$  - E $\\
   	$ L_{0,1,1} \equiv$ & $ \Delta_0   $ & $ +\left( n+1\right)\Gamma $&$   $\\
   	$ L_{1,1,1} \equiv$ & $ \Delta_0   $ & $ +\left( n+1\right)\Gamma $&$  - E $.\\
   \end{tabular} 
   $$
   \noindent
   We get minimal surfaces of general type with
   \begin{align*}
   K^2 = 16n-10, p_g = 2n, q = 0,
   \end{align*}
   and $ \degree\left( \image\varphi_{\left| K_X\right| }  \right) =2n-2 $. Moreover, $ \left| K_X\right| $ is not base point free (we will prove this in the next section~\ref{construction1.3}). Therefore, we obtain the surfaces described  in the fourth row of Theorem \ref{the main theorem}. 
   
   \subsubsection{The fixed part of the canonical system}
   \label{construction1.3}
   In this section, we show that the canonical linear system $ \left| K_X\right| $ of the surfaces in the fourth row of Theorem \ref{the main theorem} has a nontrivial fixed part. Indeed, the $ \mathbb{Z}_2^3- $cover $ \xymatrix{f: X \ar[r] &Y} $ factors through $ X_2 $, where $ X_2 $ is obtained by the $ \mathbb{Z}_2^2- $cover ramifying on $ 2L_{1,1,1}, 2L_{1,0,1} $. The linear system $ \left| K_{X_2} \right|  $ is base point free. The surface $ X $ is obtained by the $ \mathbb{Z}_2- $cover ramifying on the pullback of $ D_5 = E $ and some $ A_1 $ points. So the moving part of $ \left| K_{X} \right|  $ is the pullback of $ \left| K_{X_2} \right|  $. Therefore, the fixed part of $ \left| K_{X} \right|  $ is $ \frac{1}{2}f^{*}\left( E\right)  $. More precisely, we consider the $ \mathbb{Z}_2^3- $cover as the composition of the following $ \mathbb{Z}_2- $covers   	
   $$
   \xymatrix{X \ar[0,3]^{\mathbb{Z}_2^3}_f \ar[1,1]^{2L_{1,0,0}}_{f_3} \ar[3,0]_{\varphi_{\left| K_X \right| }}&&& Y\\
   	&X_2 \ar[r]^{2L_{1,0,1}}_{f_2} \ar[2,-1]^{\varphi_{\left| K_{X_2} \right| }}_{4:1} &X_1  \ar[-1,1]^{f_1}_{2L_{1,1,1}}&\\
   	&&&\\
   	\mathbb{P}^{2n-1}&&&} 	
   $$
   \noindent
   The first cover ramifies on $ D_2 + D_7 $ (i.e. $ B=2L_{1,1,1} $) and we get a surface $ X_1 $ with $ K_{X_1} \equiv f_1^{*}\left( -\Delta_0 + \left( n-2\right) \Gamma \right)  $. Moreover, $ f_1^*\left( E\right) = E_1 $ with $ E_1^2 = -2, g\left( E_1\right) = 0 $. The second cover ramifies on $ D_3 + D_6 $ (i.e. $ B=2L_{1,0,1} $). We have 
   \begin{align*}
   K_{X_2} \equiv f_2^{*}f_1^{*}\left( \Delta_0 + n\Gamma -E\right).
   \end{align*}
   So $ \left| K_{X_2} \right|  $ is base point free. Moreover, $ f_2^{*}\left( E_1 \right) = E_2  $ with $ E_2^2 = -4, g\left( E_2\right) = 1 $.
   The last cover ramifies on $ f_2^*f_1^{*}\left( E\right)  $ and $ 8n+6 $ nodes (i.e. $ B=2L_{1,0,0} $ ). These nodes come from the intersection points between $ D_2 $ and $D_7$, and $ D_3 $ and $D_6 $. And we obtain $ f_3^{*}\left( E_2 \right) = 2E_3  $ with $ E_3^2 = -2, g\left( E_3\right) = 1 $. In addition, by the projection formula (see \cite[Corollary 2.3]{MR3217634}), we get
   \begin{align} \label{reason 1 of non trivial fixed part of surfaces 1}
   h^0\left( K_{X} \right) = h^0\left( f_3^{*}\left( K_{X_2}\right) \right) = 2n. 
   \end{align}
   \noindent
   On the other hand, $ K_{X} \equiv f_3^{*}\left( K_{X_2}\right) + R $, where $ R $ is the ramification of $ f_3 $. Hence, 
   \begin{align}\label{reason 2 of non trivial fixed part of surfaces 1}
   K_{X} \equiv f_3^{*}\left( K_{X_2}\right) + E_3.
   \end{align}
   From $ \left( \ref{reason 1 of non trivial fixed part of surfaces 1} \right)  $ and $ \left( \ref{reason 2 of non trivial fixed part of surfaces 1} \right)  $, the elliptic curve $ E_3 $ is the fixed part of $ \left| K_{X}\right| $.
   
   \subsection{Construction 2 }
   In this section, we construct the surfaces in the last five rows of Theorem \ref{the main theorem}. 
   \subsubsection{Construction and computation of invariants}
   Let $ D_3 = \Gamma, D_4 \in \left| \Delta_0 + \Gamma \right| + \Delta_0, D_7 = \left( 2n+1\right) \Gamma $ be in $ \mathbb{F}_1$ and $ D_5, D_6 \in \left| 2\Delta_0 + 2\Gamma \right| $ be smooth curves in general position in $ \mathbb{F}_1. $ Let $ \xymatrix{f: X \ar[r] & \mathbb{F}_1} $ be a $ \mathbb{Z}_2^3- $ cover with the following branch locus
   \begin{align*}
   B= D_1 + D_2 +D_3 + D_4 + D_5 + D_6 + D_7,
   \end{align*}
   \noindent
   where $ D_1 = D_2 = 0 $. By Proposition \ref{Construction of cover}, $ L_{1,0,0} \equiv 3\Delta_0 +\left( n+3\right) \Gamma  $ and $ L_{\chi} $ is equivalent to either $ 2\Delta_0 +2\Gamma $, $ \Delta_0 +\left( n+2\right) \Gamma$ or $  \Delta_0 +\left( n+1\right) \Gamma $ for all $ L_{\chi} \neq L_{1,0,0} $. Since each $ D_{\sigma} $ is smooth and $ B $ is a normal crossings divisor, $ X $ is smooth. Furthermore, by Proposition \ref{invariants}, we get
   \begin{align*}
   2K_X &\equiv f^*\left( 2\Delta_0 + \left( 2n+1\right) \Gamma \right).
   \end{align*}
   \noindent
   This implies that $ X $ is a minimal surface of general type. Moreover, by Proposition \ref{invariants}, the invariants of $ X $ are as follows:
   \begin{align}\label{chern class 1 of surfaces 2}
   K_X^2= 16n  
   \end{align}
   \begin{align} \label{pg of surfaces 2}
   p_g\left( X\right) &= h^0\left( \Delta_0 + n\Gamma\right) =2n+1
   \end{align}
   \begin{align}\label{chi of surfaces 2}
   \chi\left( \mathcal{O}_X\right) &=2n+2.
   \end{align}
   \noindent
   From $ \left( \ref{pg of surfaces 2}\right)  $ and $ \left( \ref{chi of surfaces 2}\right)  $, we get $ q\left( X\right) = 0. $ \\
   
   We show that $ \left| K_X\right| $ is not composed with a pencil by considering the following double cover
   \begin{align*}
   \xymatrix{g_1: Y_{1} \ar[r]& \mathbb{F}_1}
   \end{align*}
   ramifying on $ D_4 + D_5 + D_6 + D_7 $. We have 
   \begin{align*}
   K_{Y_{1}}& \equiv g_1^*\left( \Delta_0 + n\Gamma \right).  
   \end{align*}
   Because $ \left| \Delta_0 + n\Gamma \right|  $ is not composed with a pencil, $ \left| K_{Y_{1}}\right|  $ is not composed with a pencil, either. This yields that $ \left| K_X\right| $ is not composed with a pencil and the degree of the canonical map is $ 8 $. 
   
   \subsubsection{The fixed part of the canonical system}
   In this section, we show that the canonical linear system $ \left| K_X\right| $ has a nontrivial fixed part. In fact, the $ \mathbb{Z}_2^3- $cover $ \xymatrix{f: X \ar[r] &Y} $ factors through $ X_2 $, where $ X_2 $ is obtained by the $ \mathbb{Z}_2^2- $cover ramifying on $ 2L_{1,1,1}, 2L_{0,1,1} $. The linear system $ \left| K_{X_2} \right|  $ is base point free. The surface $ X $ is obtained by the $ \mathbb{Z}_2- $cover ramifying on the pullback of $ D_3 = \Gamma $ and some $ A_1 $ points. So the moving part of $ \left| K_{X} \right|  $ is the pullback of $ \left| K_{X_2} \right|  $. Therefore, the fixed part of $ \left| K_{X} \right|  $ is $ \frac{1}{2}f^{*}\left( \Gamma\right)  $. More precisely, we consider the $ \mathbb{Z}_2^3- $cover as the compositions of the following $ \mathbb{Z}_2- $covers 	
   $$
   \xymatrix{X \ar[0,3]^{\mathbb{Z}_2^3}_f \ar[1,1]^{2L_{0,1,0}}_{f_3} \ar[3,0]_{\varphi_{\left| K_X \right| }}&&& Y\\
   	&X_2 \ar[r]^{2L_{0,1,1}}_{f_2} \ar[2,-1]^{\varphi_{\left| K_{X_2} \right| }}_{4:1} &X_1  \ar[-1,1]^{f_1}_{2L_{1,1,1}}&\\
   	&&&\\
   	\mathbb{P}^{2n}&&&} 	
   $$
   \noindent
   The first cover ramifies on $ D_4 + D_7 $ (i.e. $ B=2L_{1,1,1} $). We get a surface $ X_1 $ with $ K_{X_1} \equiv f_1^{*}\left( -\Delta_0 +\left( n-2\right)  \Gamma\right)$. Furthermore,  $ f_1^*\left( D_3\right) =\Gamma_1 $ with $ g\left( \Gamma_1\right) = 0 $. The second cover ramifies on $ D_5 + D_6 $ (i.e. $ B=2L_{0,1,1} $). We get surface of general type $ X_2 $ with
   \begin{align*}
   K_{X_2} \equiv f_2^{*}f_1^{*}\left( \Delta_0 +n\Gamma\right).
   \end{align*}
   Hence, $ \left| K_{X_2} \right|  $ is base point free and $ \degree\left( \image\varphi_{\left| K_{X_2}\right| }  \right) =2n-1 $. Furthermore,  $ f_2^*\left( \Gamma_1\right) =\Gamma_2 $ with $ g\left( \Gamma_2\right) = 3 $. The last cover ramifies on $ f_2^{*}f_1^{*}\left( D_3\right)  $ and $ 8n +12$ nodes (i.e. $ B=2L_{0,1,0} $ ). These nodes come from the intersection points between $ D_4 $ and $D_7$, and $ D_5 $ and $D_6 $. And we get $ f_3^*\left( \Gamma_2\right) =2\Gamma_3 $ with $ g\left( \Gamma_3\right) = 3 $. In addition, by the projection formula, we get
   \begin{align} \label{reason 1 of non trivial fixed part of surfaces 2}
   h^0\left( K_{X} \right) = h^0\left( f_3^{*}\left( K_{X_2}\right) \right) = 2n+1. 
   \end{align}
   \noindent
   On the other hand, $ K_{X} \equiv f_3^{*}\left( K_{X_2}\right) + R $, where $ R $ is the ramification of $ f_3 $. Hence, 
   \begin{align}\label{reason 2 of non trivial fixed part of surfaces 2}
   K_{X} \equiv f_3^{*}\left( K_{X_2}\right) + \Gamma_3.
   \end{align}
   Therefore, from $ \left( \ref{reason 1 of non trivial fixed part of surfaces 2} \right)  $ and $ \left( \ref{reason 2 of non trivial fixed part of surfaces 2} \right)  $, the curve $ \Gamma_3 $ is the fixed part of $ \left| K_{X}\right| $. 
   
   \subsubsection{Variations}
   By Proposition $ \ref{Construction of cover},  $ the branch locus can be imposed a point $ \left( 0, 0,0,1,1,1,1\right)$. In fact, let $ P $ be a point in $ \mathbb{F}_1 $ such that $ D_4, D_5, D_6, D_7 $ contain $ P $ with multiplicity $ 1,1,1,1 $, respectively. Let $ Y $ be the blow up of $ \mathbb{F}_1 $ at $ P $ and $ E $ be the exceptional divisor. If we abuse notation and denote $ D_3,D_4, D_5, D_6, D_7, \Delta_0, \Gamma $ their pullbacks to $ Y $, then $ D_3 = \Gamma , D_4 = 2\Delta_0 + \Gamma - E, D_5 = 2\Delta_0 + 2\Gamma - E, D_6 = 2\Delta_0 + 2\Gamma  - E  $ and $ D_7 = \left( 2n+1\right)\Gamma  - E $. Let $ \xymatrix{f: X \ar[r] & Y}  $ be a $ \mathbb{Z}^3_2- $cover with the following branch locus
   \begin{align*}
   B =D_1 + D_2 + D_3 + D_4 + D_5+ D_6 + D_7,
   \end{align*}
   \noindent
   where $ D_1 = D_2 = 0 $.  The building data is as follows:
   $$	
   \begin{tabular}{l r r r}
   	$ L_{1,0,0} \equiv$ & $ 3\Delta_0 $ & $ +\left( n+3\right) \Gamma $&$  - 2E $\\
   	$ L_{0,1,0} \equiv$ & $\Delta_0   $ &$ +\left( n+2\right)\Gamma$&$ - E $\\
   	$ L_{0,0,1} \equiv$ & $\Delta_0   $ &$ +\left( n+2\right)\Gamma$&$ - E $\\
   	$ L_{1,1,0} \equiv$ & $ 2\Delta_0   $ & $ +2\Gamma $&$ - E  $\\
   	$ L_{1,0,1} \equiv$ & $ 2\Delta_0   $ & $ +2\Gamma $&$  - E $\\
   	$ L_{0,1,1} \equiv$ & $ 2\Delta_0   $ & $ +2\Gamma $&$ - E  $\\
   	$ L_{1,1,1} \equiv$ & $ \Delta_0   $ & $ +\left( n+1\right)\Gamma $&$  - E $.\\
   \end{tabular} 
   $$
   \noindent	
   Similarly to the above, we get minimal surfaces of general type with
   \begin{align*}
   K^2 = 16n-8, p_g = 2n, q = 0, d=8,
   \end{align*}
   and $ \degree\left( \image\varphi_{\left| K_X\right| }  \right) =2n-2 $. Moreover, $  \frac{1}{2}f^{*}\left(\Gamma\right)  $ is the fixed part of $ \left| K_{X}\right|  $ and the following diagram commutes 
   $$
   \xymatrix{X \ar[0,3]^{\mathbb{Z}_2^3}_f \ar[1,1]^{2L_{0,1,0}}_{f_3} \ar[3,0]_{\varphi_{\left| K_X \right| }}&&& Y\\
   	&X_2 \ar[r]^{2L_{0,1,1}}_{f_2} \ar[2,-1]^{\varphi_{\left| K_{X_2} \right| }}_{4:1} &X_1  \ar[-1,1]^{f_1}_{2L_{1,1,1}}&\\
   	&&&\\
   	\mathbb{P}^{2n-1}&&&} 	
   $$
   \noindent
   So we obtain the surfaces in the sixth row of Theorem \ref{the main theorem}.\\
   
   Analogously, by Proposition $ \ref{Construction of cover},  $ we can put a point $ \left( 0,0,0,0,2,2,0\right)  $ into the original branch locus. In fact, let $ P $ be a point in $ \mathbb{F}_1 $ such that $ D_5, D_6$ contain $ P $ with multiplicity $ 2,2 $, respectively. Let $ Y $ be the blow up of $ \mathbb{F}_1 $ at $ P $ and $ E $ be the exceptional divisor. If we abuse notation and denote $ D_3,D_4, D_5, D_6, D_7, \Delta_0, \Gamma $ their pullbacks to $ Y $, then $ D_3 = \Gamma , D_4 = 2\Delta_0 + \Gamma , D_5 = 2\Delta_0 + 2\Gamma - 2E, D_6 = 2\Delta_0 + 2\Gamma  - 2E  $ and $ D_7 = \left( 2n+1\right)\Gamma $. Let $ \xymatrix{f: X \ar[r] & Y}  $ be a $ \mathbb{Z}^3_2- $cover with the following branch locus
   \begin{align*}
   B =D_1 + D_2 + D_3 + D_4 + D_5+ D_6 + D_7,
   \end{align*}
   \noindent
   where $ D_1 = D_2 = 0 $.  The building data is as follows:
   $$	
   \begin{tabular}{l r r r}
   	$ L_{1,0,0} \equiv$ & $ 3\Delta_0 $ & $ +\left( n+3\right) \Gamma $&$  - 2E $\\
   	$ L_{0,1,0} \equiv$ & $\Delta_0   $ &$ +\left( n+2\right)\Gamma$&$ - E $\\
   	$ L_{0,0,1} \equiv$ & $\Delta_0   $ &$ +\left( n+2\right)\Gamma$&$ - E $\\
   	$ L_{1,1,0} \equiv$ & $ 2\Delta_0   $ & $ +2\Gamma $&$ - E  $\\
   	$ L_{1,0,1} \equiv$ & $ 2\Delta_0   $ & $ +2\Gamma $&$  - E $\\
   	$ L_{0,1,1} \equiv$ & $ 2\Delta_0   $ & $ +2\Gamma $&$ - 2E  $\\
   	$ L_{1,1,1} \equiv$ & $ \Delta_0   $ & $ +\left( n+1\right)\Gamma $.&$  $\\
   \end{tabular} 
   $$
   \noindent	
   Similarly to the above, we get minimal surfaces of general type with	
   \begin{align*}
   K^2 = 16n-8, p_g = 2n, q = 1, d=8,
   \end{align*}
   and $ \degree\left( \image\varphi_{\left| K_X\right| }  \right) =2n-2 $. Furthermore, $  \frac{1}{2}f^{*}\left(\Gamma\right)  $ is the fixed part of $ \left| K_{X}\right|  $ and the following diagram commutes 
   $$
   \xymatrix{X \ar[0,3]^{\mathbb{Z}_2^3}_f \ar[1,1]^{2L_{0,1,0}}_{f_3} \ar[3,0]_{\varphi_{\left| K_X \right| }}&&& Y\\
   	&X_2 \ar[r]^{2L_{0,1,1}}_{f_2} \ar[2,-1]^{\varphi_{\left| K_{X_2} \right| }}_{4:1} &X_1  \ar[-1,1]^{f_1}_{2L_{1,1,1}}&\\
   	&&&\\
   	\mathbb{P}^{2n-1}&&&} 	
   $$
   \noindent
   Thus, we obtain the surfaces in the seventh row of Theorem \ref{the main theorem}. The Albanese pencil of these surfaces $ \xymatrix{X \ar[r]& Alb\left( X\right) } $ is the pullback of the Albanese pencil of the intermediate surface $ Z $, where $ Z $ is obtained by the $ \mathbb{Z}_2 -$cover ramifying on $ 2L_{0,1,1} $.\\
   
   Similarly, by Proposition $ \ref{Construction of cover},  $ a new branch locus can be formed by adding a point $ \left( 0,0,-1,1,2,0,1\right)  $, where $ -1 $ in the third component means the exceptional divisor $ E $ is added to $ D_3 $. In fact, let $ P $ be a point in $ \mathbb{F}_1 $ such that $ D_4, D_5, D_7$ contain $ P $ with multiplicity $ 1,2,1 $, respectively. Let $ Y $ be the blow up of $ \mathbb{F}_1 $ at $ P $ and $ E $ be the exceptional divisor. If we abuse notation and denote $ D_4, D_5, D_6, D_7, \Delta_0, \Gamma $ their pullbacks to $ Y $, then $ D_4 = 2\Delta_0 + \Gamma -E, D_5 = 2\Delta_0 + 2\Gamma - 2E, D_6 = 2\Delta_0 + 2\Gamma   $ and $ D_7 = \left( 2n+1\right)\Gamma -E$. Let $ \xymatrix{f: X \ar[r] & Y}  $ be a $ \mathbb{Z}^3_2- $cover with the following branch locus
   \begin{align*}
   B =D_1 + D_2 + D_3 + D_4 + D_5+ D_6 + D_7,
   \end{align*}
   \noindent
   where $ D_1 = D_2 = 0 $ and $ D_3 = \Gamma +E $.  The building data is as follows:
   $$	
   \begin{tabular}{l r r r}
   	$ L_{1,0,0} \equiv$ & $ 3\Delta_0 $ & $ +\left( n+3\right) \Gamma $&$  - 2E $\\
   	$ L_{0,1,0} \equiv$ & $\Delta_0   $ &$ +\left( n+2\right)\Gamma$&$  $\\
   	$ L_{0,0,1} \equiv$ & $\Delta_0   $ &$ +\left( n+2\right)\Gamma$&$ - E $\\
   	$ L_{1,1,0} \equiv$ & $ 2\Delta_0   $ & $ +2\Gamma $&$ - E  $\\
   	$ L_{1,0,1} \equiv$ & $ 2\Delta_0   $ & $ +2\Gamma $&$   $\\
   	$ L_{0,1,1} \equiv$ & $ 2\Delta_0   $ & $ +2\Gamma $&$ - E  $\\
   	$ L_{1,1,1} \equiv$ & $ \Delta_0   $ & $ +\left( n+1\right)\Gamma $&$ - E $.\\
   \end{tabular} 
   $$
   \noindent	
   Similarly to the above, we get minimal surfaces of general type with
   \begin{align*}
   K^2 = 16n -2, p_g = 2n, q = 0, d=8,
   \end{align*}
   and $ \degree\left( \image\varphi_{\left| K_X\right| }  \right) =2n-2 $. Moreover, $  \frac{1}{2}f^{*}\left(\Gamma + E\right)  $ is the fixed part of $ \left| K_{X}\right|  $ and the following diagram commutes 
   $$
   \xymatrix{X \ar[0,3]^{\mathbb{Z}_2^3}_f \ar[1,1]^{2L_{0,1,0}}_{f_3} \ar[3,0]_{\varphi_{\left| K_X \right| }}&&& Y\\
   	&X_2 \ar[r]^{2L_{0,1,1}}_{f_2} \ar[2,-1]^{\varphi_{\left| K_{X_2} \right| }}_{4:1} &X_1  \ar[-1,1]^{f_1}_{2L_{1,1,1}}&\\
   	&&&\\
   	\mathbb{P}^{2n-1}&&&} 	
   $$
   \noindent
   Therefore, we obtain the surfaces in the eighth row of Theorem \ref{the main theorem}. \\
   
   Finally, for $ n \ge 3 $ by Proposition $ \ref{Construction of cover},  $ a point $ P = \left( 0,0,-1,1,2,2,1\right)  $ can be added to the original branch locus, where $ -1 $ in the third component means the exceptional divisor is added to $ D_3 $. In fact, let $ P $ be a point in $ \mathbb{F}_1 $ such that $ D_4, D_5, D_6, D_7$ contain $ P $ with multiplicity $ 1,2,2,1 $, respectively. Let $ Y $ be the blow up of $ \mathbb{F}_1 $ at $ P $ and $ E $ be the exceptional divisor. If we abuse notation and denote $ D_4, D_5, D_6, D_7, \Delta_0, \Gamma $ their pullbacks to $ Y $, then $ D_4 = 2\Delta_0 + \Gamma -E, D_5 = 2\Delta_0 + 2\Gamma - 2E, D_6 = 2\Delta_0 + 2\Gamma - 2E  $ and $ D_7 = \left( 2n+1\right)\Gamma -E$. Let $ \xymatrix{f: X \ar[r] & Y}  $ be a $ \mathbb{Z}^3_2- $cover with the following branch locus
   \begin{align*}
   B =D_1 + D_2 + D_3 + D_4 + D_5+ D_6 + D_7,
   \end{align*}
   \noindent
   where $ D_1 = D_2 = 0 $ and $ D_3 = \Gamma +E $.  The building data is as follows:
   $$	
   \begin{tabular}{l r r r}
   	$ L_{1,0,0} \equiv$ & $ 3\Delta_0 $ & $ +\left( n+3\right) \Gamma $&$  - 3E $\\
   	$ L_{0,1,0} \equiv$ & $\Delta_0   $ &$ +\left( n+2\right)\Gamma$&$ - E $\\
   	$ L_{0,0,1} \equiv$ & $\Delta_0   $ &$ +\left( n+2\right)\Gamma$&$ - E $\\
   	$ L_{1,1,0} \equiv$ & $ 2\Delta_0   $ & $ +2\Gamma $&$ - E  $\\
   	$ L_{1,0,1} \equiv$ & $ 2\Delta_0   $ & $ +2\Gamma $&$ - E  $\\
   	$ L_{0,1,1} \equiv$ & $ 2\Delta_0   $ & $ +2\Gamma $&$ - 2E  $\\
   	$ L_{1,1,1} \equiv$ & $ \Delta_0   $ & $ +\left( n+1\right)\Gamma $&$ - E $.\\
   \end{tabular} 
   $$	
   After contracting the $ -1 $ curve arising from the fiber passing throught $ P $, we get minimal surfaces of general type with
   \begin{align*}
   K^2 = 16n-16, p_g = 2n-2, q = 1, d=8,
   \end{align*}
   and $ \degree\left( \image\varphi_{\left| K_{X}  \right|} \right) =2n-4$. 
   
   Furthermore, $  \frac{1}{2}f^{*}\left(\Gamma + E\right)  $ is the fixed part of $ \left| K_{X}\right|  $ and the following diagram commutes 
   $$
   \xymatrix{X \ar[0,3]^{\mathbb{Z}_2^3}_f \ar[1,1]^{2L_{0,1,0}}_{f_3} \ar[3,0]_{\varphi_{\left| K_X \right| }}&&& Y\\
   	&X_2 \ar[r]^{2L_{0,1,1}}_{f_2} \ar[2,-1]^{\varphi_{\left| K_{X_2} \right| }}_{4:1} &X_1  \ar[-1,1]^{f_1}_{2L_{1,1,1}}&\\
   	&&&\\
   	\mathbb{P}^{2n-3}&&&} 	
   $$
   Thus, taking $ m = n -1 $, $ m \ge 2 $, we obtain the surfaces in the last row of Theorem \ref{the main theorem}. The Albanese pencil of these surfaces $ \xymatrix{X \ar[r]& Alb\left( X\right) } $ is the pullback of the Albanese pencil of the intermediate surface $ Z $, where $ Z $ is obtained by the $ \mathbb{Z}_2 -$cover ramifying on $ 2L_{0,1,1} $.

\section*{Acknowledgments}
The author is deeply indebted to Margarida Mendes Lopes for all her help. Many thanks are also due to the anonymous referee for his/her suggestions. The author is supported by Funda\c{c}\~{a}o para a Ci\^{e}ncia e Tecnologia (FCT), Portugal under the framework of the program Lisbon Mathematics PhD (LisMath).


\Addresses

\end{document}